# EASILY RETRIEVABLE OBJECTS AMONG THE NEO POPULATION

## D. García Yárnoz,[*] J. P. Sánchez,[†] and C. R. McInnes[†]

**ABSTRACT**

Asteroids and comets are of strategic importance for science in an effort to understand the formation, evolution and composition of the Solar System. Near-Earth Objects (NEOs) are of particular interest because of their accessibility from Earth, but also because of their speculated wealth of material resources. The exploitation of these resources has long been discussed as a means to lower the cost of future space endeavours. In this paper, we consider the currently known NEO population and define a family of so-called Easily Retrievable Objects (EROs), objects that can be transported from accessible heliocentric orbits into the Earth's neighbourhood at affordable costs. The asteroid retrieval transfers are sought from the continuum of low energy transfers enabled by the dynamics of invariant manifolds; specifically, the retrieval transfers target planar, vertical Lyapunov and halo orbit families associated with the collinear equilibrium points of the Sun-Earth Circular Restricted Three Body problem. The judicious use of these dynamical features provides the best opportunity to find extremely low energy Earth transfers for asteroid material. A catalogue of asteroid retrieval candidates is then presented. Despite the highly incomplete census of very small asteroids, the ERO catalogue can already be populated with 12 different objects retrievable with less than 500 m/s of $\Delta v$. Moreover, the approach proposed represents a robust search and ranking methodology for future retrieval candidates that can be automatically applied to the growing survey of NEOs.

***Keywords:*** *Asteroids dynamics; Asteroid capture; Near-Earth Objects; Libration points*

[*] Advanced Space Concepts Laboratory, Dept. of Mechanical and Aerospace Engineering, University of Strathclyde, Glasgow G1 1XQ, UK. E-mail: daniel.garcia-yarnoz@strath.ac.uk Telephone: +44 (0)141 444 8316. Fax: +44 (0)141 552 5105
[†] Advanced Space Concepts Laboratory. Dept. of Mechanical and Aerospace Engineering, University of Strathclyde

# 1  Introduction

Recently, significant interest has been devoted to the understanding of minor bodies of the Solar System, including near-Earth and main belt asteroids and comets. NASA, ESA and JAXA have conceived a series of missions to obtain data from such bodies, having in mind that their characterisation not only provides a deeper insight into the formation of the Solar System, but also represents a technological challenge for space exploration. Near Earth Objects in particular have also stepped into prominence because of two important issues: they are among the easiest bodies to reach from the Earth and they may represent a potential impact threat.

NEOs had traditionally been classified into three families according to their orbital elements: Atens, Apollos and Amors, with Atens and Apollos being Earth-crossers, and Amors having orbits completely outside the orbit of the Earth. In recent literature (Michel, Zappala et al. 2000; Greenstreet, Ngo et al. 2011), further emphasis has been placed on the description of asteroids inside the Earth orbit, and a fourth group, the symmetric equivalent of Amors, has been added to the list. The new family has been named Atira after the first confirmed object of its kind in 2003, 163693 Atira. This is a useful classification for NEOs into 4 distinct families, and it is possible to draw some conclusions from it regarding the origin and evolution of these objects and their detectability. However, it provides little information in terms of the accessibility of their orbits.

Because of current interest in the science and exploration of NEOs, other classifications have arisen. Some of them have somewhat arbitrary or not so precise definitions: Arjunas have been defined as NEOs in extremely Earth-like orbits (Bombardelli, Urrutxuay et al. 2012), with low eccentricity, low inclination and a semi-major axis close to that of the Earth; while Brasser and Wiegert (2008) proposed a similar Small-Earth Approachers (SEA) definition for objects with diameter less than 50 m and a semi-major axis, eccentricity and incination within the ranges of [0.95 AU, 1.05 AU] , [0, 0.1] and [0°, 10°] respectively.

Other definitions concern objects that follow very particular trajectories, such as objects in horseshoe orbits, Earth's trojans, or objects that for a short period of time naturally become weakly captured by Earth, referred to as Natural Earth Satellites (NESs), or Temporarily Captured Orbiters (Granvik, Vaubaillon et al. 2011). The number of known NEOs in each of these categories is however small.

In order to provide a systematic classification of accessible objects, NASA began publishing in 2012 the Near-Earth Object Human Space Flight Accessible Target Study (NHATS) list (Abell, Barbee et al. 2012), which will be continuously updated and identifies potential candidate objects for human missions to asteroids. NEOs in NASA's NHATS list are ranked according to the number of feasible return trajectories to that object found by an automated search within certain constraints. This provides an objective, quantifiable and ordered classification of the objects in NEO space that allow feasible return missions.

Further classification involving impact hazard by NEOs have also resulted in the generation of an objective scale, the Palermo scale (Chesley, Chodas et al. 2002), for the ranking of a subset of these objects: the Potentially Hazardous Objects.

Inspired by this, and considering the growing interest in the capture of small NEOs (Sanchez and McInnes 2011; Brophy, Culick et al. 2012; Hasnain, Lamb et al. 2012), we put forward a new objective, quantifiable and ordered classification of NEOs that can be captured under certain conditions: the sub-category of Easily Retrievable Objects (EROs). EROs are defined as objects that can be gravitationally captured in bound periodic orbits around the collinear libration points $L_1$ and $L_2$ of the Sun-Earth system under a certain $\Delta v$ threshold, arbitrarily selected for this work at 500 m/s. These objects can then be ranked according to the required $\Delta v$ cost.

## 1.1 Motivation and Background

As witnesses of the early Solar System, NEOs could cast some light into the unresolved questions about the formation of planets from the pre-solar nebula, and perhaps settle debates on the origin of water on Earth or panspermian theories, among others. This scientific importance has translated into an increasing number of robotic probes sent to NEOs, and many more planned for the near future. Their low gravity well have also identified them as the only "planetary" surface that can be visited by crewed missions under NASA's flexible path plan (Augustine, Austin et al. 2009), asteroids have also become one of the feasible "planetary" surfaces to be visited by crewed missions. Science however is not the only interest of these objects and mission concepts exploring synergies with science, planetary protection and space resources utilization have started to be uttered. Examples of this are recent NASA and ESA studies on a kinetic impact demonstration mission on a binary object, DART and AIM (Murdoch, Abell et al. 2012).

Proposed technologies and methods for the deflection of potentially Earth-impacting objects have experienced significant advances, along with increasing knowledge of the asteroid population. While initially devised to mitigate the hazard posed by global impact threats, the current impact risk is largely posed by the population of small undiscovered objects (Shapiro, A'Hearn et al. 2010), and thus methods have been proposed to provide subtle changes to the orbits of small objects, as opposed to large-scale interventions such as the use of nuclear devices (Kleiman 1968). This latter batch of deflection methods, such as the low thrust tugboat (Scheeres and Schweickart 2004), gravity tractor (Edward and Stanley 2005) or small kinetic impactor (Sanchez and Colombo 2013) are moreover based on currently proven space technologies. They can therefore render the apparently ambitious scenario of manipulating asteroid trajectories a likely option for the near future.

On the other hand, the in-situ utilisation of resources in space has long been suggested as the means of lowering the cost of space missions, by means of, for example, providing bulk mass for radiation shielding or manufacturing propellant for interplanetary transfers (Lewis 1996). The development of technologies for in-situ resource utilisation (ISRU) could become a potentially disruptive innovation for space exploration and utilisation and, for example, enable large-scale space ventures that could today be considered far-fetched, such as large space solar power satellites or sustaining communities in space.

Although the concept of asteroid mining dates back to the early rocketry pioneers (Tsiolkovsky 1903), evidences of a renewed interest in the topic can be found in the growing body of literature (Baoyin, Chen et al. 2010; Sanchez and McInnes 2011; Hasnain, Lamb et al. 2012), as well as in high profile private enterprise ventures such as by Planetary Resource

s Inc[3].

With regards to the accessibility of asteroid resources, recent work by Sanchez and McInnes (2011; 2012) demonstrates that a substantial quantity of resources can indeed be accessed at relatively low energy; on the order of $10^{14}$ kg of material could potentially be harvested at an energy cost lower than that required to access resources from the surface of the Moon. More importantly, asteroid resources could be accessed across a wide spectrum of energies, and thus, current technologies could be adapted to return to the Earth's neighbourhood small objects from 2 to 30 meters diameter for scientific exploration and resource utilisation purposes.

Therefore, advances in both asteroid deflection technologies and dynamical system theory, which allow new and cheaper means of space transportation, are now enabling radically new mission concepts, such as low-energy asteroid retrieval missions (Brophy, Gershman et al. 2011). These envisage a spacecraft reaching a suitable object, coupling itself to the surface and returning it, or a portion of it, to the Earth's orbital neighbourhood. Moving an entire asteroid into an orbit in the vicinity of Earth entails obvious engineering challenges, but may also allow much more flexible resource extraction in the Earth's neighbourhood, in addition to other advantages such as enhanced scientific return.

## 2  Low Energy Transport conduits

Current interplanetary spacecraft have masses on the order of $10^3$ kg, while an asteroid of 10 meters diameter will most likely have a mass of the order of $10^6$ kg. Hence, already moving such a small object, or an even larger one, with the same ease that a scientific payload is transported would demand propulsion systems orders of magnitudes more powerful and efficient; or alternatively, orbital transfers orders of magnitude less demanding than those to reach other bodies in the solar system.

---

[3] http://www.planetaryresources.com/

Solar system transport phenomena, such as the rapid orbital transitions experienced by comets Oterma and Gehrels 3, from heliocentric orbits with periapsis outside Jupiter's orbit to apoapsis within Jupiter's orbit, or the Kirkwood gaps in the main asteroid belt, are some manifestations of the sensitivities of multi-body dynamics. The same underlying principles that enable these phenomena allow also excellent opportunities to design surprisingly low energy transfers.

It has for some time been known that the hyperbolic invariant manifold structures associated with periodic orbits around the $L_1$ and $L_2$ collinear points of the Three Body Problem provide a general mechanism that controls the aforementioned solar system transport phenomena (Belbruno and Marsden 1997; Lo and Ross 1999; Koon, Lo et al. 2000). In this paper, we seek to benefit from these mathematical constructs in order to find *low-cost* trajectories to retrieve asteroid material to the Earth's vicinity.

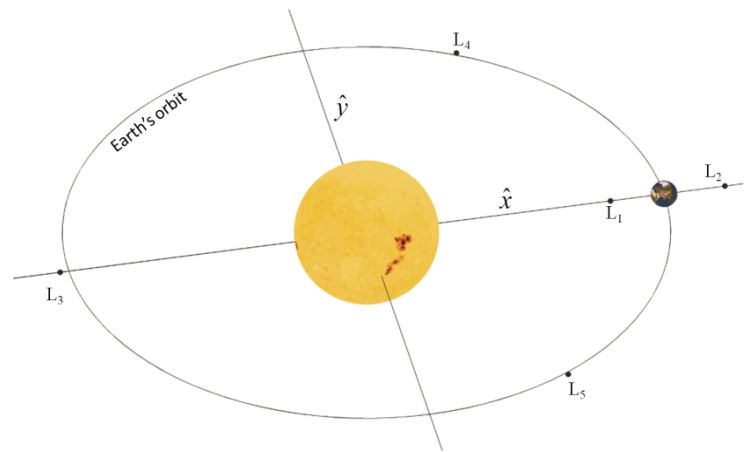

**Fig. 1**: Schematic of the CR3BP and its equilibrium points.

## 2.1 Periodic Orbits and Manifold Structure

In particular, we are interested in the dynamics concerning the Sun-Earth $L_1$ and $L_2$ points (see Fig. 1), as they are the *gate keepers* for potential ballistic capture of asteroids in the Earth's vicinity. The work in this paper assumes the motion of the spacecraft and asteroid under the gravitational influence of the Sun and Earth, within the framework of the Circular Restricted Three Body Problem (CR3BP), following closely the approach by Koon (2008). The well known equilibrium points of the system are shown in Fig. 1. The mass parameter $\mu$ considered in the paper is 3.0032080443x10$^{-6}$, which neglects the mass of the Moon. Note that the usual normalised units are used when citing Jacobi constant values.

There has been a long and intense effort to catalogue all bounded motion near the libration points of the Circular Restricted Three Body Problem (Howell 2001). The principal families of bounded motion that have been discovered are planar and vertical families of Lyapunov periodic orbits, quasi-periodic Lissajous orbits, and periodic and quasi-periodic halo orbits (Gómez, Llibre et al. 2000). Some other families of periodic orbits can be found by exploring bifurcations in the aforementioned main families (Howell 2001).

Theoretically, an asteroid transported into one of these orbits would remain near the libration point for an indefinite time. In practice, however, these orbits are unstable, and an infinitesimal deviation from the periodic orbit will make the asteroid depart asymptotically from the libration point regions. Nevertheless, small correction manoeuvres can be assumed to be able to keep the asteroid in the vicinity of the periodic orbit (Simó, Gómez et al. 1987; Howell and Pernicka 1993).

The linear behaviour of the motion near the libration points is of the type *centre* x *centre* x *saddle*, which is also a characteristic of all bounded motion near these points (Szebehely 1967). This particular dynamical behaviour ensures that, inherent to any bounded trajectory near the libration points, an infinite number of trajectories exist that asymptotically approach, or depart from, the bounded motion. Each set of trajectories asymptotically approaching, or departing, a periodic or quasi-periodic orbit near the $L_1$ or $L_2$ points forms a hyperbolic invariant manifold structure.

There are two classes of invariant manifolds: the central invariant and the hyperbolic invariant. The central invariant manifold is composed of periodic and quasi-periodic orbits near the libration points, while the hyperbolic invariant manifold consists of a stable and an unstable set of trajectories associated to the central invariant manifold. The unstable manifold is formed by the infinite set of trajectories that exponentially leaves a periodic or quasi-periodic orbit belonging to the central invariant manifold to which they are associated. The stable manifold, on the other hand, consists of an infinite number of trajectories exponentially approaching the periodic or quasi-periodic orbit.

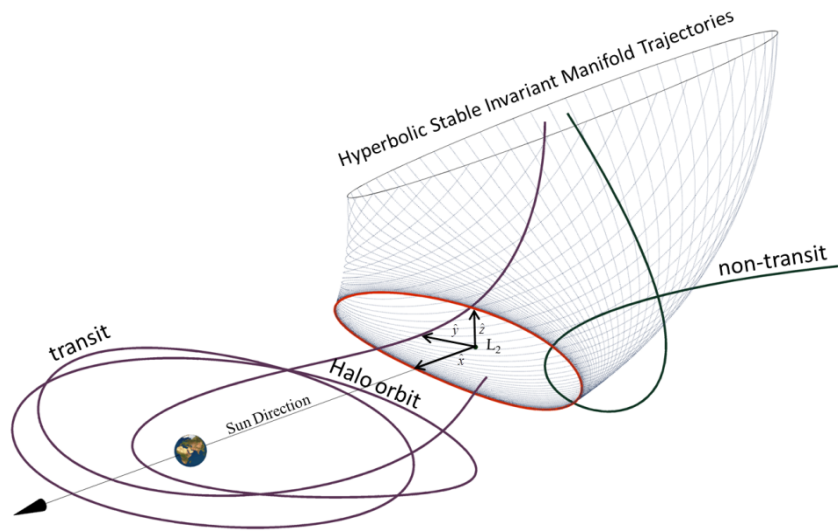

**Fig. 2:** Schematic representation of the four categories of motion near the $L_2$ point (represented by the set of axes in the figure): periodic motion around $L_2$ (i.e., halo orbit), hyperbolic invariant manifold structure (i.e., set of stable hyperbolic invariant manifold trajectories), transit trajectory and non-transit trajectory.

It is well known that the phase space near the equilibrium regions can be divided into four broad classes of motion; bound motion near the equilibrium position (i.e., periodic and quasi-periodic orbits), asymptotic trajectories that approach or depart from the latter, transit trajectories, and, non-transit trajectories (see Fig. 2). A transit orbit is a trajectory such that its motion undergoes a rapid transition between such regions. In the Sun-Earth case depicted in Fig. 2, for example, the transit trajectory approaches Earth following a heliocentric trajectory, transits through the bottle neck delimited by the halo orbit and becomes temporarily captured at Earth. An important observation from dynamical system theory is that the hyperbolic invariant manifold structure defined by the set of asymptotic trajectories forms a phase space separatrix between transit and non-transit orbits.

It follows from the four categories of motion near the libration points that periodic orbits near the Sun-Earth $L_1$ and $L_2$ points cannot only be targeted as the final destination of asteroid retrieval missions, but also as natural gateways of low energy trajectories to Earth-centred temporarily captured trajectories or transfers to other locations of the cis-lunar space, such as the Earth-Moon Lagrangian points (Lo and Ross 2001; Canalias and Masdemont 2006).

In this paper, we will focus on three distinct classes of periodic motion near the Sun-Earth $L_1$ and $L_2$ points; Planar and Vertical Lyapunov and Halo Orbits, from now on referred to as a whole as libration point orbits (LPO).

### 2.1.1 Lyapunov Orbits

As noted, the linear behaviour of the motion near the $L_1$ and $L_2$ points is of the type *centre* x *centre* x *saddle*. The *centre* x *centre* part generates a 4-dimensional central invariant manifold around each collinear equilibrium point when all energy levels are considered. In a given energy level the central invariant manifold is a 3-dimensional set of periodic and quasi-periodic solutions lying on an invariant tori, together with some stochastic regions in between (Gómez and Mondelo 2001). There exist families of periodic orbits with frequencies related to both centers: $\omega_p$ and $\omega_v$ (Alessi 2010). They are known as *planar Lyapunov family* and *vertical Lyapunov family*, see Fig. 3, and their existence is ensured by the Lyapunov centre theorem. Halo orbits are 3-dimensional periodic orbits that emerge from the first bifurcation of the planar Lyapunov family.

To generate the entire family of planar and vertical Lyapunov periodic orbits, we start by generating an approximate solution in a very close neighbourhood of the libration point (Howell 2001). This initial solution is corrected in the non-linear dynamics by means of a differential correction algorithm (Koon, Lo et al. 2008) over a suitable plane section that takes advantage of the known symmetries of these orbits (Zagouras and Markellos 1977). Once one periodic solution has been computed, the complete family can be generated by means of numerical continuation process that uses the previous solution as initial guess for a periodic orbit on which one of the dimension on the phase space has been perturbed slightly. By properly choosing the phase space direction to extend the solution by a continuation method, and by repeating the process iteratively one can build a family of periodic orbits with increasing Jacobi constant, as shown in Fig. 3.

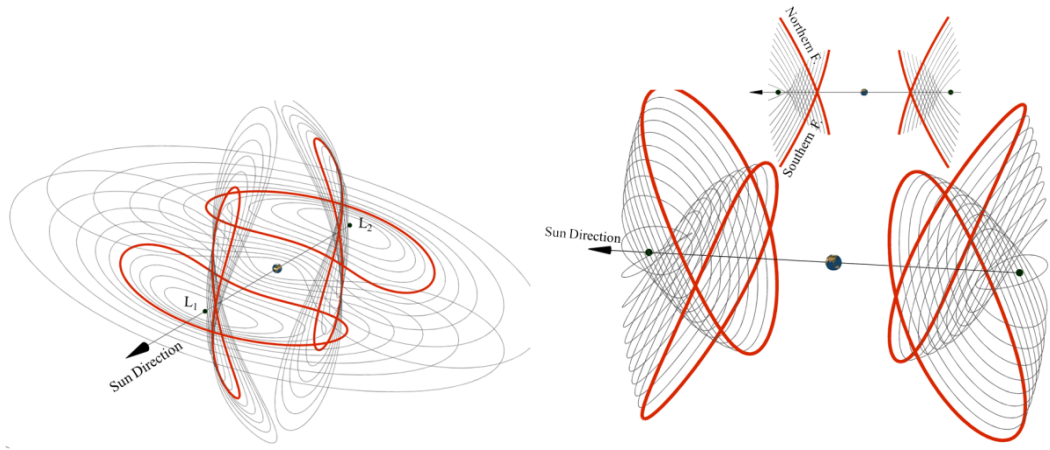

**Fig. 3:** Series of Planar and Vertical Lyapunov orbits (left) and northern and southern halo orbits (right) associated with the Sun-Earth $L_1$ and $L_2$ points. Lyapunov orbits are plotted ranging from Jacobi constant 3.0007982727 to 3.0000030032. Halo orbits are plotted ranging from Jacobi constant of 3.0008189806 to 3.0004448196. The thicker red line corresponds to a Jacobi constant of 3.0004448196, which corresponds to half the distance between the energy at equilibrium in $L_2$ and $L_3$.

### 2.1.2 Halo Orbits

The term halo orbit was coined by Robert Farquhar, who advocated the use of these orbits near the Earth-Moon $L_2$ point to obtain a continuous communication relay with the far side of the Moon during the Apollo programme (Farquhar 1967).

As previously noted, this type of orbit emerges from a bifurcation in the planar Lyapunov orbits. As the amplitude of planar Lyapunov orbits increase, eventually a critical amplitude is reached where the planar orbits become vertical critical, as defined by Hénon (1973), and new three-dimensional families of periodic orbits bifurcate. Thus, the minimum possible size for Halo orbits in the Sun-Earth system is approximately (240 x 660) x $10^3$ km at $L_1$ and (250 x 675) x $10^3$ km at $L_2$, sizes denoting the maximum excursion from the libration point in the $x$ and $y$ directions respectively. At the bifurcation point, two symmetric families of halo orbits emerge at each libration point, here referred to as the northern and southern family depending on whether the maximum $z$ displacement is achieved in the northern (i.e., $z > 0$) or southern (i.e., $z < 0$) direction, respectively (see Fig. 3).

Similarly to planar and vertical Lyapunov, the set of halo orbits, also shown in Fig. 3, was computed by means of the continuation of a predictor-corrector process. The initial seed was computed by means of Richardson (1980) third order approximation of a halo orbit. A differential corrector procedure is used to trim Richardson's prediction and obtain the smallest halo possible (Zagouras and Markellos 1977; Koon, Lo et al. 2008). We then continue the process by feeding the next iteration with a prediction of a slightly larger displacement in $z$. Repeating this process provides a series of halo orbits with increasing energy, or decreasing Jacobi constant.

## 3 Asteroid Retrieval Opportunities

In the past few years, several space missions have already attempted to return samples from the asteroid population, e.g., Hayabusa (Kawaguchi, Fujiwara et al. 2008), and others are planned for the near future[4]. As shown by Sanchez and McInnes (2011; 2012), given the low transport cost expected for the most accessible objects, we could also envisage the possibility to return to Earth entire small objects with current or near-term technology. The main challenge resides on the difficulties inherent in the detection of these small objects. Thus, for example, only 1 out of every million objects with diameter between 5 to 10 meters is currently known and this ratio is unlikely to change significantly in the coming years (Veres, Jedicke et al. 2009).

---

[4] http://www.nasa.gov/topics/solarsystem/features/osiris-rex.html (last accessed 02/05/12)

In this section then, we will focus our attention on the surveyed population of asteroids in search of the most accessible candidates for near-term asteroid retrieval missions by means of invariant hyperbolic stable manifold trajectories, the so called EROs.

For this purpose, a systematic search of capture candidates among catalogued NEOs was carried out, selecting the $L_1$ and $L_2$ regions as the target destination for the captured material. This gives a grasp and better understanding of the possibilities of capturing entire NEOs or portions of them in a useful orbit, and demonstrates a method that can be applied to categorise newly discovered small bodies in the future when detection technologies improve, and rank them according to their retrievability.

### 3.1 Invariant Manifold Trajectories to $L_1$ and $L_2$

In order to provide a simple but robust method for categorizing EROs, the design of the transfer from the asteroid orbit to the $L_1$ and $L_2$ LPO consists of a ballistic arc, with two impulsive burns at the start and end, intersecting a hyperbolic stable invariant manifold asymptotically approaching the desired periodic orbits. This paper only considers the inbound leg of a full capture mission.

Planar Lyapunov, vertical Lyapunov, and Halo orbits around $L_1$ and $L_2$ generated with the methods described in the previous section were considered as target orbits. The invariant stable manifold trajectories were computed by perturbing the target orbit periodic solutions around the Lagrangian point on the stable eigenvector direction (Koon, Lo et al. 2008) by a magnitude of $10^{-6}$, in normalized units. These initial conditions were propagated backwards in the Circular Restricted 3-Body Problem until they reached the desired fixed section in the Sun-Earth rotating frame. We refer to this propagation time as the manifold transfer time. The section was arbitrarily selected as the one forming an angle of $\pm\pi/8$ with the Sun-Earth line ($\pi/8$ for the $L_2$ orbits, see Fig. 4, the symmetrical section at $-\pi/8$ for those targeting $L_1$). This corresponds roughly to a distance to Earth of the order of 0.4 AU, where the gravitational influence of the planet is considered small. No additional perturbations were considered in the backward propagation.

In this analysis, Earth is assumed to be in a circular orbit 1 AU away from the Sun. This simplification allows the conditions of the manifold trajectories (and in particular in the selected section) to be independent of the insertion time into the final orbit. The only exception is the longitude of the perihelion, i.e., the sum of the right ascension of the ascending node and the argument of perihelion, which varies with the insertion time with respect to a reference time with the following relation:

$$(\Omega + \omega) = (\Omega_{REF} + \omega_{REF}) + \frac{2\pi}{T}(t - t_{REF}) \tag{1}$$

where $\Omega_{REF}$ and $\omega_{REF}$ are the right ascension of the ascending node and the argument of perihelion at the $\pm\pi/8$ section for an insertion into a target orbit at reference time $t_{REF}$, and $T$ is the period of the Earth. For orbits with non-zero inclination, the argument of perihelion of the manifolds is also independent of the insertion time and the above equation indicates a variation in $\Omega$. However, in the case of planar Lyapunov with zero inclination, $\Omega$ is not defined and an arbitrary value of zero can be selected, resulting in the equation representing a change in argument of perihelion.

The transfer between the NEO orbit and the manifold is then calculated as a heliocentric Lambert arc of a restricted two-body problem with two impulsive burns, one to depart from the NEO, the final one for insertion into the manifold, with the insertion constrained to take place before or at the $\pm\pi/8$ section.

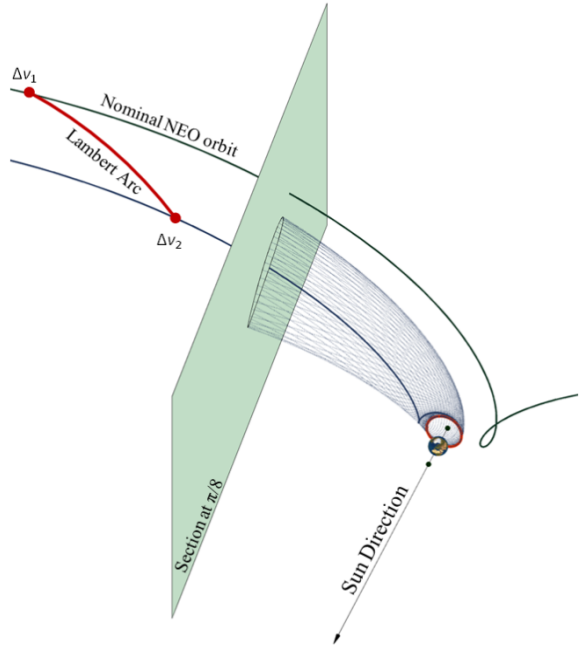

**Fig. 4:** Schematic representation of a transfer to $L_2$

Thus, the problem can be defined with 5 variables: the Lambert arc transfer time, the manifold transfer time, the insertion date at the target periodic orbit, the energy of the final orbit, and a fifth parameter determining the point in the target orbit where the insertion takes place.

The benefit of such an approach is that the asteroid is asymptotically captured into a bound orbit around a collinear Lagrangian point, with no need for a final insertion burn at arrival. All burns are performed far from Earth, so no large gravity losses need to be taken into account. Furthermore, this provides additional time for corrections, as the dynamics in the manifold are "slow" when compared to a traditional hyperbolic approach. Finally, this type of trajectory is then easily extendable to a low-thrust trajectory if the burns required are small.

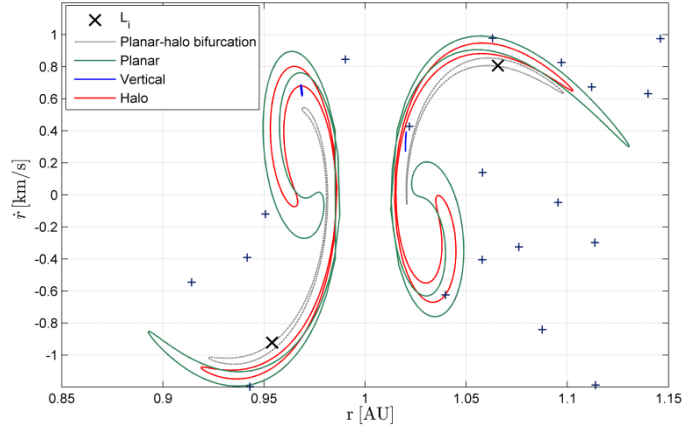

**Fig. 5:** Projection of the manifolds onto the $r-\dot{r}$ phase space for a Jacobi constant of 3.0004448196. The manifolds are represented at their intersection with a plane forming a ±π/8 angle with the Sun-Earth line in the rotating frame. Manifolds on the left correspond to $L_1$, on the right to $L_2$. Capture candidates are indicated with a + marker.

The shape of the manifolds projected onto the $r-\dot{r}$ phase space (with *r* being the radial distance from the Sun) at the intersection with the ±π/8 section is shown in Fig. 5 for a particular Jacobi constant. For an orbit with exactly the energy of $L_1$ or $L_2$, the intersection is a single point; while for lower Jacobi constants, the shape of the intersection is a closed loop. The intersection corresponding to the bifurcation between planar and halo orbits is also plotted. A few capture candidate asteroids have been included in the plot (+ markers) at the time of their intersection with the π/8 plane around their next closest approach to Earth. It is worth noting that the epoch of the next encounter, and thus of the intersection, is different for each particular asteroid. In a planar case, this would already provide a good measure of the distance of the asteroid to the manifolds. However, as we are considering the 3D problem, information on the *z* component or the inclination would also be necessary.

Figure 6 provides a more useful representation of the manifolds in terms of perihelion, aphelion radius and inclination for the two collinear points. The point of bifurcation between the planar Lyapunov and halo orbits, when they start growing in inclination, can easily be identified. Halo orbits extend a smaller range in aphelion and perihelion radius when compared to planar Lyapunov orbits. Vertical Lyapunov orbits have even smaller excursions in radius from a central point, as can already be seen in the smaller loops of vertical Lyapunov orbits in Fig. 5, but they extend to much lower values of the Jacobi constant and cover a wider range of inclinations.

Several asteroids are also plotted with small markers in the graphs. Their Jacobi constant $J$ is approximated by the Tisserand parameter as defined in Eq. (2).

$$J \approx \frac{1}{a} + 2\sqrt{a(1-e^2)}\cos i \qquad (2)$$

where $a$, $e$ and $i$ are the semi-major axis (in AU), eccentricity and inclination of the asteroid orbit.

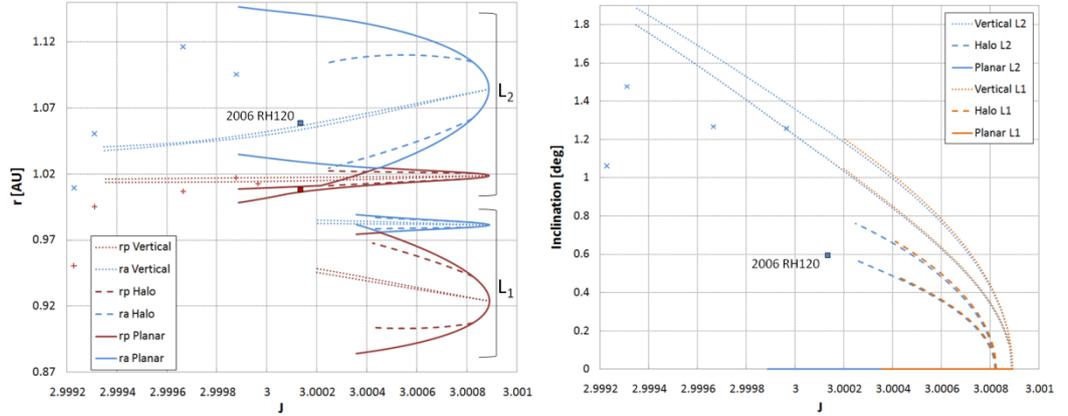

**Fig. 6:** Minimum and maximum perihelion and aphelion radius (left) and inclination (right) of the manifolds leading to planar Lyapunov, vertical Lyapunov and halo orbits around $L_1$ and $L_2$.

This illustrates the proximity to the manifolds of a number of NEOs. In particular, asteroid 2006 RH120 has been highlighted, due to its proximity to the $L_2$ manifolds. From these graphs and ignoring any phasing issues, it can already be identified as a good retrieval mission candidate, as its perihelion and aphelion radius are close to or within the range of all three types of considered manifolds, and its inclination also lies close to the halo orbit manifolds. The manifold orbital elements appear to be a good filter to prune the list of NEOs to be captured.

## 3.2 Asteroid Catalogue Pruning

For the calculation of capture opportunities, the NEO sample used for the analysis is JPL's Small Bodies Database[5], downloaded as of 27th of July of 2012. This database represents the catalogued NEOs up to that date, and as such it is a biased population, most importantly in size, as already noted. A large number of asteroids of the most ideal size for capture have not yet been detected, as current detection methods favour larger asteroids. Secondly, there is an additional detection bias related to the type of orbits, with preference for Amors and Apollos in detriment to Atens, as objects in Aten orbits spend more time in the exclusion zone due to the Sun.

Even with this reduced list, it is a computationally expensive problem and preliminary pruning becomes necessary. Previous work by Sanchez et al. (2012) showed that the number of known asteroids that could be captured from a hyperbolic approach with a total Δ$v$ less than 400 m/s is of the order of 10. Although the hyperbolic capture approach in their work and the manifold capture is inherently different, the number of bodies that could be captured in manifold orbits at low cost is expected to be of the same order. Without loss of generality, it is possible to immediately discard NEOs with semi-major axis (and thus energy) far from the Earth's, as well as NEOs in highly inclined orbits. However, a more systematic filter needed to be devised.

As a first approximation of the expected total cost in terms of Δ$v$, a bi-impulsive cost prediction with both burns assumed at aphelion and perihelion was implemented. Either of the two burns is also responsible for correcting the inclination. The Δ$v$ required to modify the semi-major axis can be expressed as:

$$\Delta v_a = \sqrt{\mu_s \left( \frac{2}{r} - \frac{1}{a_f} \right)} - \sqrt{\mu_s \left( \frac{2}{r} - \frac{1}{a_0} \right)} \qquad (3)$$

where $\mu_S$ is the Sun's gravitational constant, $a_0$ and $a_f$ are the initial and final semi-major axis before and after the burn, and $r$ is the distance to the Sun at which the burn is made (perihelion or aphelion distance). On the other hand the Δv required to modify the inclination at either apsis can be approximated by:

$$\Delta v_i = 2\sqrt{\frac{\mu_s}{a_0}r^*} \sin(\Delta i / 2) \qquad (4)$$

where Δi is the required inclination change, and r* corresponds to the ratio of perihelion and aphelion distance if the burn is performed at aphelion, or its inverse if performed at perihelion.

---

[5] http://ssd.jpl.nasa.gov/sbdb.cgi (last accessed 27/07/12)

Note that these formulas are only first order approximations intended for the pruning of the database, and they will not be used to calculate the final transfers. In particular, the plane change is only valid for small changes in inclination and large deviations from the values provided by the filter are expected to be observed for large inclinations. Nevertheless, we are interested in low cost transfers which imply a small plane change, so the approximation is acceptable. Also, these formulas only take into consideration the shape and inclination of the orbits, ignoring the rest of the orbital elements: right ascension of the ascending node and argument of pericentre. It is then implicitly assumed that the line of nodes coincides with the line of apsis and the inclination change can be performed at pericentre or apocentre.

The total estimated cost for pruning is then calculated as:

$$\Delta v_t = \sqrt{\Delta v_{a1}^2 + \Delta v_{i1}^2} + \sqrt{\Delta v_{a2}^2 + \Delta v_{i2}^2} \qquad (5)$$

with one burn performed at each of the apsis, and one of the two inclination change Δv assumed zero.

The estimated transfer Δv corresponds thus to the minimum of four cases: aphelion burn modifying perihelion and inclination followed by a perihelion burn modifying aphelion, perihelion burn modifying aphelion and inclination followed by an aphelion burn modifying perihelion, and the equivalent ones in which the inclination change is done in the second burn.

For simplicity, the target manifold final perihelion, aphelion and inclination values are selected as ranges or bands obtained from Fig. 6. For example, planar Lyapunov manifolds at $L_2$ correspond to a range of *[$r_p$, $r_a$, i]* ∈ *[1.00-1.02, 1.02-1.15,0]*, or *[1.01-1.02, 1.025-1.11,0.59-0.78]* for halo manifolds at $L_2$. Note that the inclination range for halos was given as the one that corresponds to the highest energy. This is due to the fact that most candidate asteroids have higher energies than the manifolds, and the lowest cost is assumed to take place where the energy difference is minimum. In the case of vertical Lyapunov orbits, due to the narrow ranges and strong dependency with *J*, polynomial fits for *[$r_p$, $r_a$, i]* as a function of *J* were used.

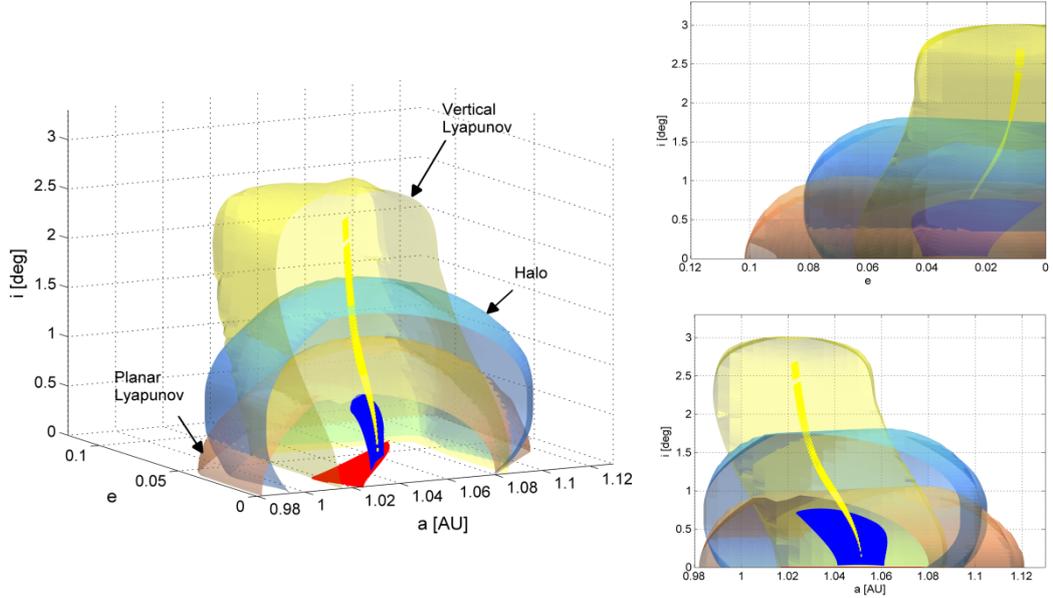

**Fig. 7:** Regions in the orbital element space with total estimated cost for capture into an LPO around $L_2$ below 500 m/s. The manifolds corresponding to the LPOs are plotted in solid colours.

With this filter, it is then possible to calculate the regions of a three-dimensional orbital element space (in semi-major axis, eccentricity and inclination) than can potentially be captured under a certain Δ*v* threshold. These regions are plotted in Fig. 7 for transfers to LPOs around $L_2$ with a Δ*v* of 500 m/s, and any asteroid with orbital elements inside them could in principle be captured at that cost. The figure shows a three-dimensional view of the surfaces that delimit the regions for planar Lyapunov, vertical Lyapunov and halo, as well as two-dimensional projections in the *a-i* and *e-i* planes. There is a significant overlap between the regions of different LPO target orbits; therefore, it is expected that several asteroids would allow low-cost captures to more than one family of LPO. A similar plot can be generated for the case of $L_1$. Figure 8 presents the regions for $L_1$ and $L_2$ compared to the definitions of the 4 families of NEOs. Objects from all four families seem to be adequate candidates for the new category of Easily Retrievable Objects, particularly the ones closed to the Apollo-Amor and Aten-Atira divides. The boundaries for the Small-Earth Approachers subset is also depicted with a dashed line, and shows that this definition is not particularly useful for the purpose of pruning candidates for asteroid retieval.

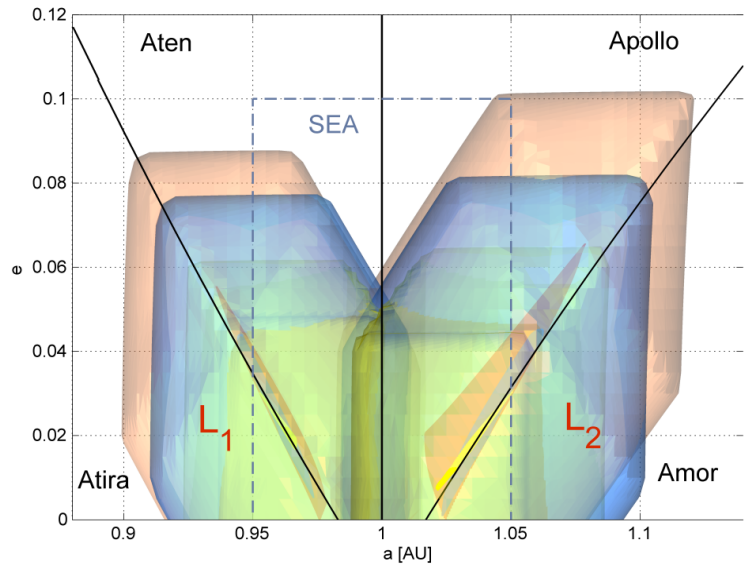

**Fig. 8:** Semi-major axis and eccentricity map of the capturable regions for $L_1$ and $L_2$. The boundaries of the main 4 families of NEO objects and the Small-Earth Approachers subset are also indicated. The manifold orbital elements are enclosed in the capturable regions and closely follow the Apollo-Amor and Aten-Atira divides.

The filter approximation provides in general a lower bound $\Delta v$ estimate, as it ignores any phasing issues, and assumes the burns can be performed at apocentre or pericentre. Moreover, there is no guarantee, and in fact it is quite unlikely, that a combination of the extremes of the ranges of $[r_p, r_a, i]$ used in the filter correspond to proper manifold trajectories. Finally, the plane change does not include a modification in right ascension of the ascending node. Although the final $\Omega$ can be tuned by modifying the phasing with the Earth, this is not completely free as the final insertion will take place around a natural close approach of the asteroid with the planet. The combination of this constrained phasing and the plane change will also incur in additional costs. North and south halo obits provide two opportunities with opposite $\Omega$ for each transfer, which should result in two different costs, while the filter provides a single value.

For a few cases, with high initial inclination and associated plane change cost, the filter can over-estimate the $\Delta v$. As the inclination increases, solutions splitting the large plane change into the two burns can potentially result in a lower cost. In cases where the filter favours solutions with larger burns at pericentre, it can also incur in higher costs estimation for the plane change than the optimal solution.

## 4   Capture Transfer Results

As the main objective is to catalogue objects that can be captured under a threshold of 500 m/s, we will focus on the filtered asteroids with estimated $\Delta v$ below 1 km/s as provided by the filter, to be on the conservative side. For each of these NEOs, feasible capture transfers with arrival dates in the interval 2016-2100 were obtained. The NEO orbital elements are only considered valid until their next close encounter with Earth. The Lambert transfers between the asteroid initial orbit and the manifolds were optimised using EPIC, a global optimisation method that uses a stochastic search blended with an automatic solution space decomposition technique (Vasile and Locatelli 2009). Single objective optimisations with total transfer $\Delta v$ as the cost function were carried out. Trajectories obtained with EPIC were then locally optimised with MATLAB's built-in constrained optimisation function fmincon.  Lambert arcs with up to 3 complete revolutions before insertion into the manifold were considered. For cases with at least one complete revolution, the two possible solutions of the Lambert problem were optimised. This implies that 7 full problem optimisations needed to be run for each NEO. In order to limit the total duration of the transfers, the insertion into the manifold was arbitrarily constrained to take place not earlier than 1000 days before the $\pm\pi/8$ section during the global search. This constraint was released in the local optimisation.

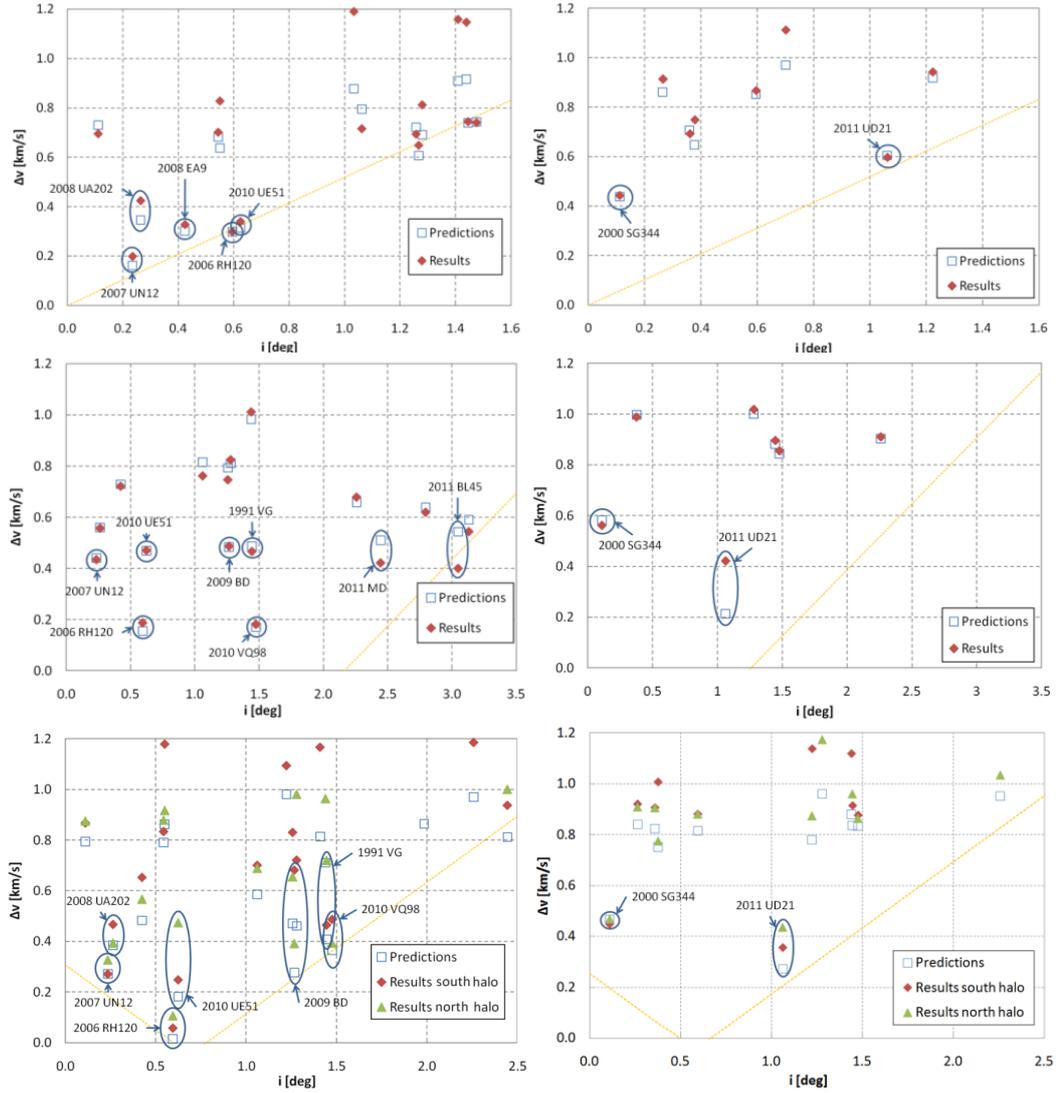

**Fig. 9:** Filter cost estimates and results of the optimisation for planar Lyapunov (top), vertical Lyapunov (middle) and halo orbits (bottom) around $L_2$ (left) and $L_1$ (right). Dotted lines indicate the cost of changing just the inclination.

Figure 9 plots the results of the optimisation for $L_2$ and $L_1$ together with the estimates. It can be observed that the filter provides in general a good approximation of the total cost to be expected. As expected, the larger the inclination, the larger the deviation of the results from the predicted cost by the filter. It is nevertheless a useful tool to select candidates and prioritise lists of asteroids for optimisation, and to quickly predict if any newly discovered asteroid is expected to have low capture costs. Dotted lines have been added to the plot as indicators of the cost of performing just the inclination change at a circular orbit at 1 AU. Predicted and optimised results are expected to fall above or close to these lines. EROs with capture costs smaller than 500 m/s are identified in the plots.

**Table 1:** NEO characteristics for transfer trajectories with Δ$v$ below 500 m/s. The type of transfer is indicated by a 1 or 2 indicating $L_1$ or $L_2$ plus the letter P for planar Lyapunov, V for vertical Lyapunov, and Hn or Hs for north and south halo.

| Rank # | | a [AU] | e | i (Pravec, Scheirich et al.) | MOID [AU] | Diameter [m] | Type | Δ$v$ [km/s] |
|---|---|---|---|---|---|---|---|---|
| 1 | 2006 RH120 | 1.033 | 0.024 | 0.595 | 0.0171 | 2.3- 7.4 | 2Hs | 0.058 |
| | | | | | | | 2Hn | 0.107 |
| | | | | | | | 2V | 0.187 |
| | | | | | | | 2P | 0.298 |
| 2 | 2010 VQ98 | 1.023 | 0.027 | 1.476 | 0.0048 | 4.3-13.6 | 2V | 0.181 |
| | | | | | | | 2Hn | 0.393 |
| | | | | | | | 2Hs | 0.487 |
| 3 | 2007 UN12 | 1.054 | 0.060 | 0.235 | 0.0011 | 3.4-10.6 | 2P | 0.199 |
| | | | | | | | 2Hs | 0.271 |
| | | | | | | | 2Hn | 0.327 |
| | | | | | | | 2V | 0.434 |
| 4 | 2010 UE51 | 1.055 | 0.060 | 0.624 | 0.0084 | 4.1-12.9 | 2Hs | 0.249 |
| | | | | | | | 2P | 0.340 |
| | | | | | | | 2V | 0.470 |
| | | | | | | | 2Hn | 0.474 |
| 5 | 2008 EA9 | 1.059 | 0.080 | 0.424 | 0.0014 | 5.6-16.9 | 2P | 0.328 |
| 6 | 2011 UD21 | 0.980 | 0.030 | 1.062 | 0.0043 | 3.8-12.0 | 1Hs | 0.356 |
| | | | | | | | 1V | 0.421 |
| | | | | | | | 1Hn | 0.436 |
| 7 | 2009 BD | 1.062 | 0.052 | 1.267 | 0.0053 | 4.2-13.4 | 2Hn | 0.392 |
| | | | | | | | 2V | 0.487 |
| 8 | 2008 UA202 | 1.033 | 0.069 | 0.264 | 2.5·10$^{-4}$ | 2.4- 7.7 | 2Hn | 0.393 |
| | | | | | | | 2P | 0.425 |
| | | | | | | | 2Hs | 0.467 |
| 9 | 2011 BL45 | 1.033 | 0.069 | 3.049 | 0.0040 | 6.9-22.0 | 2V | 0.400 |
| 10 | 2011 MD | 1.056 | 0.037 | 2.446 | 0.0018 | 4.6-14.4 | 2V | 0.422 |
| 11 | 2000 SG344 | 0.978 | 0.067 | 0.111 | 8.3·10$^{-4}$ | 20.7-65.5 | 1P | 0.443 |
| | | | | | | | 1Hs | 0.449 |
| | | | | | | | 1Hn | 0.468 |
| 12 | 1991 VG | 1.027 | 0.049 | 1.445 | 0.0037 | 3.9-12.5 | 2Hs | 0.465 |
| | | | | | | | 2V | 0.466 |

Table 1 shows the EROs with capture costs lower than the selected Δ$v$ threshold. Twelve asteroids of the whole NEO catalogue can be retrieved at this cost, ten of them around $L_2$ plus two Atens around $L_1$. The table provides the orbital elements, minimum orbit intersection distance according to the JPL Small Bodies Database, and an estimate of the size of the object. This estimate is calculated with the following relation (Chesley, Chodas et al. 2002):

$$D = 1329\,km \times 10^{-H/5}\, p_v^{-1/2} \qquad (6)$$

where the absolute magnitude $H$ is provided in the JPL database, and the albedo $p_v$ is assumed to range from 0.05 (dark) to 0.50 (very bright icy object).

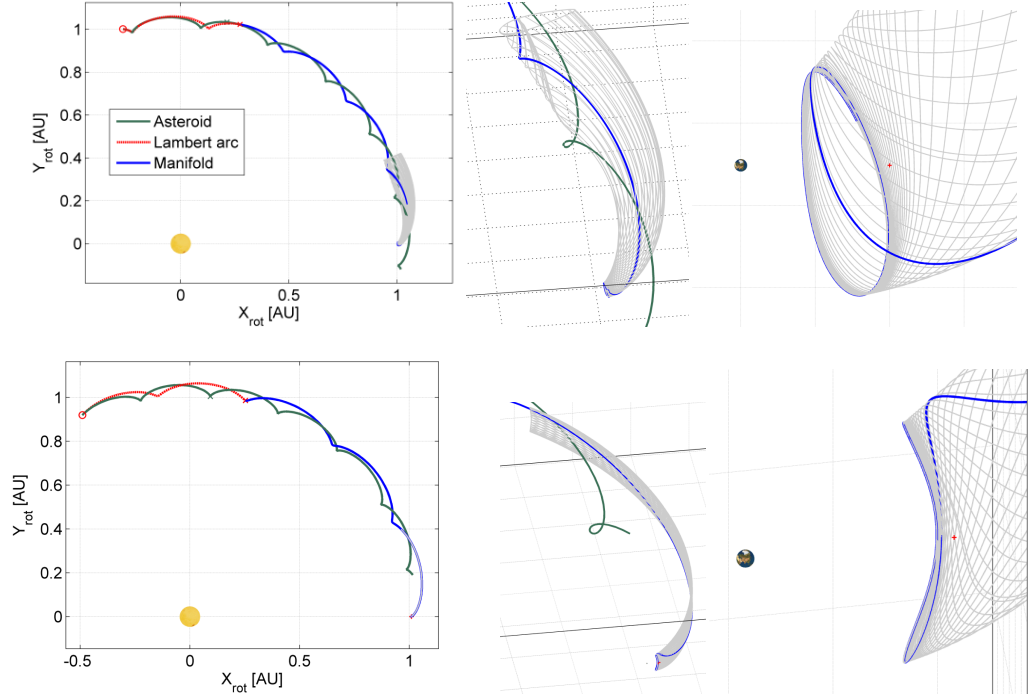

**Fig. 10:** Capture trajectories for asteroid 2006 RH120 to a south halo (top) and vertical Lyapunov (bottom). The unperturbed original orbit of the asteroid is plotted in dark green. Sun and Earth are not to scale; they are plotted 10 times their size.

As expected, planar Lyapunov orbits are optimal for lower inclination NEOs, while NEOs with higher inclination favour transfers to vertical Lyapunov. Figure 10 shows two example trajectories in a co-rotating frame where the Sun-Earth line is fixed for a transfer of asteroid 2006 RH120 to LPOs around $L_2$. Both trajectories correspond to the same close approach of the asteroid to Earth in 2028. Close-ups of the final parts of the trajectory are plotted in a three-dimensional view in order to appreciate the shape of the final orbit and manifolds.

Multiple trajectories were found for each asteroid, lasting between 2.2 and 10 years. Table 2 presents the best trajectory for each type of target orbit for $L_2$ and $L_1$. The cheapest transfer, below 60 m/s, corresponds to a trajectory inserting asteroid 2006

RH120 into a halo orbit. Solutions to planar and vertical Lyapunov orbits were also found for 2006 RH120 at higher costs. This agrees well with the interpretation of Fig. 6. The pruning method was also predicting that this transfer would be the cheapest, with a minimum estimated $\Delta v$ of 15 m/s. It is important to emphasise that the total $\Delta v$ comprises both burns at departure from the asteroid and insertion into the manifold. The NEO orbit may intersect the manifold directly, and in that case the transfer to the target orbit can be done with a single burn, as in this particular asteroid.

The total duration of the transfers range from 3 to 7.5 years. For the longer transfers it is possible to find faster solutions with less revolutions in the Lambert arc at a small $\Delta v$ penalty.

**Table 2:** Capture trajectories and mass estimates for the best trajectory of each type.

|  |  | Date [yy/mm/dd] | | | J manifold | Total Durat. [yr] | $\Delta v$ [m/s] | | $I_{sp}$ = 300s | |
|---|---|---|---|---|---|---|---|---|---|---|
|  |  | Asteroid departure | Manifold insertion | $L_i$ arrival |  |  | Dep | Ins | Mass [ton] | Ø [m] |
| 2006 RH120 | 2Hs | 21/02/01 | 21/02/01 | 28/08/05 | 3.000421 | 7.51 | 58 | 0 | 153.6 | 4.83 |
| 2006 RH120 | 2Hn | 23/05/11 | 24/02/20 | 28/08/31 | 3.000548 | 5.31 | 52 | 55 | 82.3 | 3.92 |
| 2010 VQ98 | 2V | 35/02/14 | 35/09/01 | 39/11/15 | 3.000016 | 4.75 | 177 | 4 | 46.8 | 3.25 |
| 2007 UN12 | 2P | 13/10/22 | 13/10/22 | 21/02/19 | 3.000069 | 7.33 | 199 | 0 | 42.3 | 3.14 |
| 2011 UD21 | 1Hs | 37/11/20 | 38/07/03 | 42/07/19 | 3.000411 | 4.66 | 149 | 207 | 21.9 | 2.52 |
| 2011 UD21 | 1V | 36/07/20 | 38/11/16 | 41/06/21 | 3.000667 | 4.92 | 226 | 196 | 17.9 | 2.36 |
| 2011 UD21 | 1Hn | 39/10/24 | 40/06/15 | 43/08/30 | 3.000504 | 3.85 | 210 | 226 | 17.2 | 2.33 |
| 2000 SG344 | 1P | 24/02/11 | 25/03/11 | 27/06/18 | 3.000357 | 3.35 | 195 | 248 | 16.8 | 2.04 |

# 5 Discussion

## 5.1 Overview of the catalogue of EROs

All identified EROs are of small size (perhaps with the exception of 2000 SG344), which is ideal for a technology demonstrator retrieval mission. In fact, seven of them fit the SEA definition by Brasser and Wiegert (2008). They showed, focusing on object 1991 VG, that the orbit evolution of these type of objects is dominated by close encounters with Earth, with a chaotic variation in the semi-major axis over long periods of time. A direct consequence of this is that reliable capture transfers can only be designed with accuracy over one synodic period, before the next encounter with Earth changes the orbital elements significantly. The fact that EROs are close to the hyperbolic manifolds makes them a particularly interesting subset of NEOs with regards to dynamics, since they represent objects with potential for high sensitivity of gravitational perturbation during these future Earth encounters. One could argue that finely tuning these encounters could also be used to shepherd these objects into trajectories that have a lower cost to be inserted into a manifold (Sanchez and McInnes 2011).

The NEOs in Table 1 are well-known, and there has been speculation about the origin of a few of them, including the possibility that they were man-made objects (spent upper stages) or lunar ejecta after an impact (Tancredi 1997; Chodas and Chesley 2001; Brasser and Wiegert 2008; Kwiatkowski, Kryszczynska et al. 2009). In particular object 2006 RH120 has been thoroughtly studied (Kwiatkowski, Kryszczynska et al. 2009; Granvik, Vaubaillon et al. 2011), as it was a temporarily captured object that was considered the "second moon of the Earth" until it finally escaped the Earth in July 2007. Granvik shows that the orbital elements of 2006 RH120 changed from being an asteroid of the Atens family pre-capture, to an Apollo post-capture, having followed what we refer in this paper to as a transit orbit inside Earth's Hill sphere, and thus its must have orbited inside the separatrix surface of the hyperbolic stable manifold. An additional object in the list, 2007 UN12, is also pointed out by Granvik as a possible candidate to become a TCO.

Regarding their accessibility, a recent series of papers (Adamo, Giorgini et al. 2010; Barbee, Espositoy et al. 2010; Hopkins, Dissel et al. 2010) considered up to 7 of the above objects as possible destinations for the first manned mission to a NEO (and the other 5 were not discovered at the time). They proposed human missions during the same close approaches as the capture opportunities calculated. However, the arrival dates at the asteroids are later than the required departure date for the capture, so their outbound legs could not apply to our proposed capture trajectories. An additional study by Landau and Strange (2011) presents crewed mission trajectories to over 50 asteroids. It shows that a mission to 6 of the considered asteroids is possible with a low-thrust $\Delta v$ budget between 1.7 and 4.3 km/s. The costs presented are for a return mission of a spacecraft with a dry mass of 36 tons (including habitat) in less than 270 days. A longer duration robotic mission with a final mass at the NEO of less than 10 tons and a manifold capture as proposed here would result in much lower fuel costs as the thrust-to-mass ratio increases. Moreover, eleven of our 12 capturable objects appear in the top 25 of NASA's NHATS list as of September 2012, seven of them in the top 10. This indicates that the objects found by our pruning and optimisation are indeed easily accessible, even if the outbound part of the trajectory was not considered in our calculation.

## 5.2 Retrievable mass limit with current space technology

The results presented in the previous section could be used to calculate a limit in the mass that can be retrieved with current space technology. In order to obtain a first estimate of the mass and size of the asteroids that can be captured, we can consider a basic system mass budget exercise. The Keck study report for asteroid retrieval (Brophy, Gershman et al. 2011) proposes a mission involving a spacecraft of 5500 kg dry mass and 8100 kg of propellant already at the NEO encounter. With a spacecraft of those characteristics, the total asteroid mass that could be transferred with the trajectories described in this paper is close to 400 tons. However, the launch mass required would be close to 16 tons. Such a high launch mass would imply either a long escape strategy from LEO, or a heavy launcher not yet developed, or multiple launches and assembly in space. We can consider a more modest mission of the size of Cassini (2442 kg dry mass and 3132 kg propellant mass[6]) at the NEO. A full system budget would require a larger fuel mass at launch to deliver the spacecraft to the target, and thus an analysis of the outbound leg. However, preliminary analysis for asteroid 2006 RH120, performed in the frame of an asteroid deflection demonstrator mission, show trajectories with low departure velocities from Earth (well below 1 km/s) and transfer $\Delta v$ budgets lower than 450 m/s. These figures translate into a spacecraft of 6300 kg departing Earth with an escape velocity of around 500 m/s, within the capabilities of current launch systems such as Ariane 5 ECA. Multiple burn escape strategies from a HEO orbit are also feasible.

---

[6] http://saturn.jpl.nasa.gov/multimedia/products/pdfs/cassini_msn_overview.pdf (last accessed 05/09/12)

Assuming the Cassini-like mass budgets, results are appended for each trajectory on Table 2 for a standard high-thrust propulsion system. The total mass for a high thrust engine of specific impulse ($I_{sp}$) 300s ranges from 17 to about 154 tons, which represents 3 to 28 times the wet mass of the spacecraft at arrival to the NEO. The trajectories presented assume impulsive burns, so in principle they are not suitable for low-thrust transfers. However, due to their low $\Delta v$ and long time of flight, transformation of these trajectories to low-thrust is in principle feasible, and will be considered in future work. If a similar cost trajectory could be flown with a low-thrust engine of higher specific impulse (e.g., 3000s) the asteroid retrieved mass would be over ten times that of the high-thrust case, up to an impressive 1500 tons or over 10 m diameter in the case of a hypothetical transfer from the orbit of 2006 RH120 to a halo orbit.

For an average NEO density of 2.6 gr/cm$^3$ (Chesley, Chodas et al. 2002), the equivalent diameter of the asteroid that can be captured is also included in the table. This shows that reasonably sized boulders of 2-5 m diameter, or entire small asteroids of that size, could be captured with this method. The capture of entire bodies of larger size is still challenging, but the derived size of a few of the candidates fall actually within this range. The ERO 2000 SG344, with a derived size in the range of 20 to 65 meters, is the only asteroid that completely fails to meet the capturable range shown in Table 2, even with the higher specific impulse.

Regarding the safety of such a project, there could be a justified concern regarding the possibility of an uncontrolled re-entry of a temporary captured asteroid into Earth atmosphere. A migration through the unstable invariant manifold leading towards the inner region around Earth could result on homoclinic or heteroclinic transits between $L_1$ and $L_2$ (Koon, Lo et al. 2000), some of which intersect the planet. An active control would be required to ensure that all deviation from the target periodic or quasi-periodic orbits is in the direction of the unstable manifolds leading to the outside (for $L_2$) or inside ($L_1$) heliocentric regions. It is however a less serious concern due to the small size of the considered EROs. Objects smaller than 5 meters have a low impact energy (specially if we consider the lower velocity impacts that would result from a transit orbit when compare to a hyperbolic trajectory), and a relatively high impact frequency with Earth (Chesley, Chodas et al. 2002). Statistically, one object of a similar size impacts the Earth every 1-3 years with limited consequences. If larger objects were considered, additional mitigation measures would be required. The Keck study report (Brophy, Culick et al. 2012) suggests a Moon orbit as the final destination for their captured object, to circumvent this problem.

## 5.3 Method Limitations

One of the first objections that can be raised to the approach presented involves some of the simplifications in the model. The main simplifying assumptions are placing the Earth in a circular orbit, assuming Keplerian propagation for the NEOs orbital elements until the next close encounter with Earth without considering any uncertainties in their ephemerides, and not including other types of perturbations, in particular the Moon third body perturbation. While the influence of the first two assumptions on the general behaviour of the trajectories should be relatively small, and the transfers obtained can be used as first guesses for a local optimisation with a more complex model with full Earth and NEOs ephemerides, not including the Moon as a perturbing body can have a much greater influence. Granvik (2011) shows that the Moon plays an important role in the capture of TCO, and the trajectories of the manifolds would be also affected by it. The lunar third body perturbation can also strongly influence the stability of LPOs, in particular large planar Lyapunov orbits, and it could render some of them unsuitable as target orbits. Trajectories calculated with full dynamics may no longer be optimal, the final orbits are no longer periodical in an Elliptical Restricted 3-Body Problem, and they can also become unstable. A control strategy would be required to maintain a captured object in an orbit around a Lagrangian point. However, the asymptotic behaviour of the manifolds and the type of NEOs that can be captured are not expected to change. The family of EROs presented are also of large scientific interest as they are the most likely candidates to suffer natural transitions through the $L_1/L_2$ regions and migrations between NEO families. Other perturbations, such as the changes in the orbit of small bodies affected by solar radiation pressure are of little importance within the timescales considered.

Even if unstable, the target libration point orbits presented can serve though as either observation points for the temporarily captured EROs, or as gateways to other Sun-Earth-Moon system orbits of interest, through the transit orbits inside Earth's Hill sphere and heteroclinic connections between libration points. Other capture possibilities, e.g. by means of a single or double lunar swingby, or multiple resonant Earth swingbys, have not been studied and are outside of the scope of this paper, but they could potentially increase the number of retrievable objects available.

## 6  Conclusions

The possibility of capturing a small NEO or a segment from a larger object would be of great scientific and technological interest in the coming decades. It is a logical stepping stone towards more ambitious scenarios of asteroid exploration and exploitation, and possibly the easiest feasible attempt for humans to modify the Solar System environment outside of Earth, or attempt a large-scale macro-engineering project.

This paper has shown that the retrieval of a full asteroid is well within today's technological capabilities, and that there exists a series of objects with potential to be temporarily captured into libration point orbits. We define these objects as Easily Retrievable Objects (EROs). These are objects whose orbits lie close to a stable hyperbolic invariant manifold such that a small $\Delta v$ transfer may link the nominal trajectory of the asteroid with an assymptotic trajectory leading to a periodic orbit near the Sun-Earth $L_1/L_2$ points. Under certain conditions, these transfers can be achieved with transfers cost below 500 m/s. Indeed, the paper presents a list of 12 EROs, with a total of 25 trajectories to periodic orbits near $L_2$ and 6 near $L_1$ below a cost of 500 m/s, and the number of these objects is expected to grow considerably in the coming years. The lowest cost is of 58 m/s to transfer asteroid 2006 RH120 to a halo southern family with a single burn on 1st February 2021. All the capture transfer opportunities to Earth's vicinity have been identified for the currently catalogued NEOs during the next 30 years, and enable capture of bodies within 2-5 meters diameter with low propellant costs.

Taking advantage of these transfer opportunities and the unique dynamical characteristics of the identified EROs, the science return of asteroid missions can be greatly improved, and the utilisation of asteroid resources may become a viable mean of providing substantial mass in Earth orbit for future space ventures. Despite the largely incomplete survey of very small objects, the current known population of asteroids provides a good starting platform to begin with the search for easily retrievable objects. With this goal, a robust methodology for systematic pruning of a NEO database and optimisation of capture trajectories through the hyperbolic invariant stable manifold into different types of LPO around $L_1$ and $L_2$ has been implemented and tested.

The proposed method can be easily automated to prune the NEO database on a regular basis, as the number of EROs in orbits of interest is expected to grow with the new efforts in asteroid detection. Any new occurrence of a low-cost candidate asteroid can be optimised to obtain the next available phasing, transfer opportunities and the optimal target LPO.

Moreover, Sun-Earth LPOs can also be considered as natural gateways to the Earth system. Thus, the problem to transfer an asteroid to an Earth or Moon centred orbit can be decoupled into the initial phase of inserting the asteroid into a stable invariant manifold and then providing the very small manoeuvres required to continue the transit into the Earth system. While a method to find optimal LPO capture trajectories and possible targets has been defined in this paper, the transit trajectories can potentially allow the asteroid to move to the Earth-Moon $L_1/L_2$ or other locations within cis-lunar space taking advantage of heteroclinic connections between collinear points.

## 7  Acknowledgements

The authors wish to acknowledge Elisa Maria Alessi for her valuable comments and inputs to this work. The work was carried out  making use of the Faculty of Engineering

High Performance Computer Facility, University of Strathclyde; and was supported by European Research Council grant 227571 (VISIONSPACE).